\documentclass[11pt,a4paper,reqno]{amsart}

\usepackage{amsmath,amsthm}
\usepackage{amssymb}

\usepackage{cases}

\theoremstyle{theorem}

\theoremstyle{definition}

\theoremstyle{remark}

\usepackage{url}
\usepackage[colorlinks,linkcolor=blue,citecolor=blue]{hyperref}

\title{More on Faulhaber's theorem for sums of powers of integers}

\author[J.L. Cereceda]{Jos\'e Luis Cereceda}
\address{%
        Collado Villalba, 28400 -- Madrid, Spain}
\email{jl.cereceda@movistar.es}

\begin{document}

\begin{abstract}
In this note we consider the theorem established in arXiv:1912.07171 concerning the sums of powers of the first $n$ positive
integers, $S_k = 1^k + 2^k + \cdots + n^k$, and show that it can be used to demonstrate the classical theorem of Faulhaber
for both cases of odd and even $k$.
\end{abstract}

\maketitle

\section{Introduction}

For integers $n \geq 1$ and $k \geq 0$, denote $S_k = \sum_{i=1}^{n} i^k$. It is well-known that $S_k$ can be expressed in the
so-called Faulhaber form (see, e.g., \cite{beardon,cere1,cere2,cere3,cere4,edwards1,edwards2,knuth,kotiah,krishna,witmer})
\begin{align}
S_{2k} & = S_2 \big[ b_{k,0} + b_{k,1} S_1 + b_{k,2} S_{1}^2 + \cdots + b_{k,k-1} S_{1}^{k-1} \big], \label{f1} \\
S_{2k+1} & =  S_{1}^2 \big[ c_{k,0} + c_{k,1} S_1 + c_{k,2} S_{1}^2 + \cdots +  c_{k,k-1} S_{1}^{k-1} \big], \label{f2}
\end{align}
where $b_{k,j}$ and $c_{k,j}$ are non-zero rational coefficients for $j =0,1,\ldots,k-1$ and $k \geq 1$. In particular, $S_3 = S_{1}^2$. We
can write \eqref{f1} and \eqref{f2} more compactly as
\begin{align*}
S_{2k} & = S_2  F_{2k}(S_1),  \\
S_{2k+1} & =  S_{1}^2 F_{2k+1}(S_1), 
\end{align*}
where both $F_{2k}(S_1)$ and $F_{2k+1}(S_1)$ are polynomials in $S_1$ of degree $k-1$. For later convenience, we also quote the
relationship between $S_{2}^2$ and $S_1$, namely
\begin{equation}\label{ship}
S_{2}^2 = \frac{1}{9} S_{1}^{2} (1 + 8 S_{1}).
\end{equation}

Recently, Miller and Trevi\~{n}o \cite{miller} derived the following formulas for $S_{2r+1}$ and $S_{2r+2}$:
\begin{align}
S_{2r+1} & = \frac{r+1}{2} \left( S_{r}^2 - \sum_{i=r}^{2r-1} d_i S_i \right),  \label{mt1}  \\
S_{2r+2} & = \frac{(r+1)(r+2)}{2r+3} \left( S_r S_{r+1} -  \sum_{i=r+1}^{2r} e_i S_i \right),  \label{mt2}
\end{align}
from which, by using mathematical induction, they concluded (see \cite[Theorem 1]{miller}) that, for $k \geq 1$, there exists a
polynomial $g_k \in \mathbb{Q}[x,y]$ such that $g_k(0,0)=0$ and $S_k = g_k(S_1,S_2)$.

In this note we argue that, actually, a slight reformulation of the said theorem enables one to demonstrate the theorem of Faulhaber
embodied in equations \eqref{f1} and \eqref{f2} above. Indeed, as will become clear in the next section, the formulas in \eqref{mt1}
and \eqref{mt2} can be used to generate recursively the Faulhaber polynomials in \eqref{f2} and \eqref{f1}, respectively. Before going
further, it should be noticed that, throughout this work, we adopt the convention of expressing $S_{2}^2$ as in equation \eqref{ship},
so that we consider that, formally, $S_{2}^2$ does {\it not\/} depend on $S_2$. According to \eqref{f1}, this implies that, if $k$ is even,
$S_{k}^2$ is a polynomial in $S_1$ (see \cite[Corollary 3.2]{beardon}). For example, we will write $S_5$ as $\frac{1}{3} (4 S_{1}^3 -
S_{1}^2 )$, and {\it not\/} as $\tfrac{1}{2}(3S_{2}^2 - S_{1}^2)$. Assuming this convention, it turns out that, for odd $k$, say $k =
2r+1$, we can make the right-hand side of \eqref{mt1} to depend exclusively on $S_1$, so that the resulting polynomial for $S_{2r+1}$
obtained from \eqref{mt1} yields the Faulhaber polynomial $S_{2r+1} = S_{1}^2 F_{2r+1}(S_1)$. On the other hand, starting with
\eqref{mt2}, one can indeed show that $S_{2r+2} = g_{2r+2}(S_1,S_2)$, but this relationship is, again, nothing more than that given
by the Faulhaber polynomial $S_{2r+2} = S_2 F_{2r+2}(S_1)$.

\section{Generators of Faulhaber polynomials}

To support our claim, we make use of the following formula, which gives us the product of the power sums $S_k$ and $S_m$
(with $k,m \geq 1$):
\begin{equation}\label{prod1}
S_k S_m = \frac{1}{k+1} \sum_{j=0}^{k/2} B_{2j} \binom{k+1}{2j} S_{k+m+1-2j} +
\frac{1}{m+1} \sum_{j=0}^{m/2} B_{2j} \binom{m+1}{2j} S_{k+m+1-2j},
\end{equation}
where the $B_j$'s are the Bernoulli numbers and where the upper summation limit $k/2$ denotes the greatest integer less than or equal
to $k$. Equation \eqref{prod1} was stated as Theorem 1 in \cite{mac}, where, incidentally, it is further observed that it was known to Lucas
by 1891. For $k =m$, \eqref{prod1} reduces to
\begin{equation}\label{square}
S_k^2 = \frac{2}{k+1} \sum_{j=0}^{k/2} B_{2j} \binom{k+1}{2j} S_{2k+1-2j}, \quad k \geq 1.
\end{equation}

From \eqref{square}, we then obtain
\begin{equation}\label{prod2}
S_{2r+1} = \frac{r+1}{2}S_{r}^2 - \sum_{j=1}^{r/2} B_{2j} \binom{r+1}{2j} S_{2r+1-2j}, \quad r \geq 1,
\end{equation}
where we can see that the summation in the right-hand side of \eqref{prod2} involves only power sums $S_i$ with odd index $i$. Thus, according
to \eqref{f2}, every $S_i$ can be put as $S_{1}^2$ times a polynomial in $S_1$. Likewise, from \eqref{f1}, \eqref{f2}, and \eqref{ship}, it turns
out that $S_{r}^2$ can always be expressed as $S_{1}^2$ times a polynomial in $S_1$ of degree $r-1$. Hence, it follows that $S_{2r+1}$
in equation \eqref{prod2} must factorize as the product of $S_{1}^2$ times a polynomial in $S_1$ of degree $r-1$, namely, the Faulhaber form
$S_{2r+1} = S_{1}^2 F_{2r+1}(S_1)$. As a simple example, we may use \eqref{prod2} to evaluate $S_9$. For $r=4$, equation \eqref{prod2}
reads
\begin{equation*}
S_9 = \frac{5}{2} S_{4}^2 + \frac{1}{6}S_5 - \frac{5}{3}S_7.
\end{equation*}
Now, since $S_4 = \frac{1}{5} (6 S_1 S_2 - S_2)$, and taking into account \eqref{ship}, we obtain
\begin{equation*}
S_{4}^2 = \frac{1}{225} S_{1}^2 - \frac{4}{225} S_{1}^3 - \frac{4}{15} S_{1}^4 + \frac{32}{25} S_{1}^5 .
\end{equation*}
Thus, noting that $S_5 = \frac{1}{3} (4 S_{1}^3 - S_{1}^2 )$ and $S_7 = \tfrac{1}{3} (6S_{1}^4 - 4S_{1}^3 + S_{1}^2$),
we finally get the Faulhaber polynomial
\begin{equation*}
S_9 = \frac{1}{5}S_{1}^2 \big( -3 + 12 S_1 - 20 S_{1}^2 + 16 S_{1}^3 \big).
\end{equation*}

On the other hand, by taking $k =r$ and $m=r+1$ in \eqref{prod1}, and solving for $S_{2r+2}$, we find that
\begin{equation}\label{prod3}  
S_{2r+2} = \frac{(r+1)(r+2)}{2r+3} \left( S_r S_{r+1} - B_{r+1} S_{r+1} - \sum_{j=1}^{r/2} h_{r,j} B_{2j} S_{2r+2-2j} \right),
\,\,\, r \geq 1,
\end{equation}
where
\begin{equation*}
h_{r,j} =  \frac{1}{r+1}\binom{r+1}{2j} + \frac{1}{r+2}\binom{r+2}{2j} .
\end{equation*}
Note that the summation in the right-hand side of \eqref{prod3} involves only sums $S_i$ with even index $i$. Furthermore, the single term
$B_{r+1}S_{r+1}$ only survives when $r+1$ is even. Regarding the product $S_r S_{r+1}$, it is obvious that one of the indices $r$ or $r+1$
is even. Therefore, invoking \eqref{f1}, we conclude that $S_{2r+2}$ in equation \eqref{prod3} must factorize as the product of $S_2$ times a
polynomial in $S_1$ of degree $r$, namely, the Faulhaber form $S_{2r+2} = S_2 F_{2r+2}(S_1)$. As another concrete example, let us evaluate
$S_{10}$ using \eqref{prod3}. For $r=4$, equation \eqref{prod3} reads
\begin{equation*}
S_{10} = \frac{30}{11} \left( S_4 S_5 + \frac{7}{60} S_6 - \frac{3}{4} S_8 \right).
\end{equation*}
Now, substituting $S_4 = \frac{1}{5} (6 S_1 S_2 - S_2)$, $S_5 = \frac{1}{3} (4 S_{1}^3 - S_{1}^2 )$ , $S_6 = \frac{1}{7}(S_2 - 6S_1 S_2
+ 12 S_{1}^2 S_2$), and $S_8 = \frac{1}{15}(-3S_2 + 18 S_1 S_2 - 40 S_{1}^2 S_2 + 40 S_{1}^3 S_2)$ into the last equation, we get the
Faulhaber polynomial
\begin{equation*}
S_{10} = \frac{1}{11} S_2 \big( 5 - 30 S_1 +68 S_{1}^2 - 80 S_{1}^3 + 48 S_{1}^4 \big).
\end{equation*}

It should then be clear that the formulas for $S_{2r+1}$ and $S_{2r+2}$ in \eqref{prod2} and \eqref{prod3} or, equivalently, the formulas in
\eqref{mt1} and \eqref{mt2} act as generators of the Faulhaber polynomials, provided that both $S_{2r+1}$ and the square $S_{2r}^2$ are
expressed in terms of $S_1$. Armed with the formulas in \eqref{prod2} and \eqref{prod3}, it is then a trivial matter to inductively prove the
Faulhaber theorem given in equations \eqref{f1} and \eqref{f2}. A proof of this kind based on equations like \eqref{prod2} and \eqref{prod3}
was given elsewhere \cite{cere1}.

For the sake of completeness, it is worth observing that the Faulhaber polynomials can also be obtained by means of the identities\footnote{
Identity \eqref{i1} appears as Theorem 2 in \cite{mac} (see also \cite[Equation (6)]{acu} and \cite[Equation (17)]{kotiah}). Regarding
identity \eqref{i2}, it can be readily obtained from \cite[Equation (22)]{kotiah} (see also \cite[Equation (4.3)]{cere2}).}
\begin{equation}\label{i1}
S_1^k = \frac{1}{2^{k-1}} \sum_{j=0}^{\frac{k-1}{2}} \binom{k}{2j+1} S_{2k-1-2j}, \quad k \geq 1,
\end{equation}
and
\begin{equation}\label{i2}
S_2 S_1^k = \frac{1}{3 \cdot 2^k} \sum_{j=0}^{\frac{k+1}{2}} \frac{2k+3-2j}{2j+1} \binom{k+1}{2j} S_{2k+2-2j}, \quad k \geq 1.
\end{equation}
From \eqref{i1} and \eqref{i2}, it follows that
\begin{equation}\label{i3}
S_{2r+1} = \frac{2^r}{r+1} S_{1}^{r+1} - \frac{1}{r+1} \sum_{j=1}^{r/2} \binom{r+1}{2j+1} S_{2r+1-2j}, 
\end{equation}
and
\begin{equation}\label{i4}
S_{2r+2} = \frac{1}{2r+3}\left( 3 S_2 (2S_{1})^{r} - \sum_{j=1}^{\frac{r+1}{2}} \frac{2r+3-2j}{2j+1}
\binom{r+1}{2j} S_{2r+2-2j} \right),
\end{equation}
respectively. Note  that the summation in the right-hand side of \eqref{i3} [\eqref{i4}] involves only odd [even] indexed power sums $S_i$,
Therefore, starting with $r =1$, one can recursively use \eqref{i3} [\eqref{i4}] to get the Faulhaber polynomials $S_{2r+1} = S_{1}^2
F_{2r+1}(S_1)$ [respectively, $S_{2r+2} = S_{2}F_{2r+2}(S_1)$].

\section{Conclusion}

In \cite[Theorem 1]{miller}, Miller and Trevi\~{n}o deduced from \eqref{mt1} and \eqref{mt2} that, for $k \geq 1$, there exist a polynomial
$g_k \in \mathbb{Q}[x,y]$ such that $g_k(0,0)=0$ and $S_k = g_k(S_1,S_2)$. The point raised in this note is that any such polynomial in $S_1$
and $S_2$ can always be reduced to the Faulhaber form in equations \eqref{f1} and \eqref{f2}. Indeed, as we have shown here, the above
polynomials \eqref{mt1} and \eqref{mt2} can be used to generate recursively the Faulhaber polynomials $S_{2r+1} = S_{1}^2 F_{2r+1}(S_1)$
and $S_{2r+2} = S_{2}F_{2r+2}(S_1)$, respectively, provided we adhere to the convention in \eqref{ship}

On the other hand, as discussed in \cite{miller}, there may be other possible ways of expressing $S_{k}$ in terms of power sums of lower degree.
In this sense, it is pertinent to recall the remarkable result achieved by Beardon (see \cite[Theorem 6.2]{beardon}), according to which, for each
pair of integers $i$ and $j$ with $1 \leq i < j$, there is a unique, non-constant irreducible polynomial $T_{ij}$ in two variables $x$ and $y$, with
integer coefficients, such that $T_{ij}(S_i, S_j) =0$. As a simple example, we have the relation (\cite[Equation (1.4)]{beardon})
\begin{equation*}
T(S_1, S_2) = 0,  \quad \text{with} \quad T(x,y) = 8x^3 + x^2 - 9y^2,
\end{equation*}
which is just relation \eqref{ship}. Moreover, it was further shown there (see \cite[Theorem 7.1]{beardon}) that the relation $T_{i,j}(S_i, S_j)$ is
separable if, and only if, $i=1$. This is, of course, in agreement with the polynomial form in \eqref{f2}. A result already anticipated by Faulhaber in
1631 \cite{edwards1}.


\vspace{.5cm}

\begin{thebibliography}{99}


\bibitem{acu} Acu, D. (1988). Some algorithms for the sums of integer powers. \textit{Math. Mag.} 61(3):189--191.


\bibitem{beardon} Beardon, A. F. (1996). Sums of powers of integers. \textit{Amer. Math. Monthly.} 103(3):201--213.


\bibitem{cere1} Cereceda, J. L. (2012). (In Spanish). Teorema de Faulhaber sobre las sumas de potencias. {\it La Gaceta de la
RSME}. 15(1):149--169.


\bibitem{cere2} Cereceda, J. L. (2013). Averaging sums of powers of integers and Faulhaber polynomials. {\it Ann. Math. Inform.}
42:105--117.

\bibitem{cere3} Cereceda, J. L. (2014). A determinant formula for sums of powers of integers. {\it Int. Math. Forum}
9(17):785--795.


\bibitem{cere4} Cereceda, J. L. (2015). Explicit form of the Faulhaber polynomials. \textit{College Math. J.} 46(5):359--363.


\bibitem{edwards1} Edwards, A. W. F. (1982). Sums of powers of integers: a little of the history. \textit{Math. Gaz.} 66(435):22--28.


\bibitem{edwards2} Edwards, A. W. F. (1986). A quick route to sums of powers. \textit{Amer. Math. Monthly.} 93(6):451--455.


\bibitem{knuth} Knuth, D. E. (1993). Johann Faulhaber and sums of powers. \textit{Math. Comp.} 61:277--294.


\bibitem{kotiah} Kotiah, T. C. T. (1993). Sums of powers of integers---A review. \textit{Int. J. Math. Educ. Sci. Technol.}
24(6):863--874.


\bibitem{krishna} Krishnapriyan, H. K. (1995). Eulerian polynomials and Faulhaber's result on sums of powers of integers.
\textit{College Math. J.} 26(2):118--123.


\bibitem{mac} MacDougall, J. A. (1988). Identities relating sums of powers of integers -- An exercise in generalization.
{\it Australian Senior Math. J.} 2(1):53--62.


\bibitem{miller} Miller, S. J. and Trevi\~{n}o, E. (2019). On the sum of $k$-th powers in terms of earlier sums. Available
at: arXiv:1912.07171v1.


\bibitem{witmer} Witmer, E. E. (1935). The sums of powers of integers. \textit{Amer. Math. Monthly.} 42(9):540--548.


\end{thebibliography}
\end{document}